\documentclass[preprint]{elsarticle}

\usepackage{amssymb,amsmath, color}
\usepackage[mathscr]{eucal}
\usepackage{eufrak} \usepackage{mathrsfs}
\usepackage{comment}

\usepackage{fancyhdr}
\makeatletter
\def\ps@pprintTitle{%
 \let\@oddhead\@empty
 \let\@evenhead\@empty
 \def\@oddfoot{}%
 \let\@evenfoot\@oddfoot}
\makeatother
\newtheorem{thm}{Theorem}  

\newtheorem{lemma}[thm]{Lemma}
\newtheorem{cor}[thm]{Corollary}

\newproof{pf}{Proof}
\newdefinition{rmk}{Remark}

\DeclareMathOperator{\sig}{sig}
\newcommand{\Cc}{\mathcal{C}}
\newcommand{\sD}{\mathscr{D}}

\begin{document}

\begin{frontmatter}

\title{The distance signatures of  the incidence graphs of affine resolvable designs }

\author[JM]{Jianmin Ma}
\ead{jianminma@yahoo.com}

\address[JM]{ College of Math \& Information Science    and Hebei Key Lab of Computational Mathematics \& Applications
 \\  Hebei Normal University, Shijiazhuang, Hebei 050016, China}


\begin{abstract}
In this note, we determined the distance signatures of the incidence matrices of affine resolvable designs. This proves a conjecture by
Kohei Yamada. 
\end{abstract}

\begin{keyword}
Distance-matrix \sep  signature \sep affine divisible design


\MSC[2008] 05C50
15A18
\end{keyword}
\date{Jan. 22, 2015}

\end{frontmatter}

\section{Introduction}

Let $A$ be a $d\times d$ symmetric matrix over the real field $\mathbb{R}$. Let $n_{+}(A)$, $n_{-}(A)$, and $n_0(A)$ be the number (counting multiplicity for repeat values) of positive, negative and zero eigenvalues of $A$, respectively.  So $d= n_{+}(A) + n_{-}(A) + n_0(A)$.  The {\em signature (inertia)} of $A$ is the triple
$(n_{+}(A) , n_{-}(A),  n_0(A))$, denoted by $\sig(A)$.

For any graph $G$, we use the same the letter $G$ to denote its vertex set. The distance   $\partial(x,y)$ of vertices $x,y\in G$ is the length of a shortest path between them. The distance matrix $D= D(G)$ is formed by indexing the rows and columns with the vertex set $G$ and defining the $(x,y)$ entry to be $\partial(x,y)$. 
Following \cite{Graham}, the signature $\sig(D)$ is called the {\em distance signature} of  graph $G$.   We also write $\sig(G)$ for $\sig(D)$ and $n_{\pm,0}(D)$ for $n_{\pm,0}(G)$.
Since $D$ is real and symmetric, it has real eigenvalues. Because the trace of $D(G)$ is zero,  $n_{+}(G)$ and $n_{-}(G)$ are bounded above by $|G|-1$.
 
In \cite{Graham},  Graham and Lov\'{a}sz  remarked that if it was not even known which graphs $G$ has $n_{-}(G) = |G|-1$ or whether there is a graph for which $n_{+}(G) > n_{-}(G)$.
{Azarija} \cite{Azarija} showed  that the family of conference graphs have $n_{+}(G) > n_{-}(G)$. If $G$ is a conference graph on  $v$ vertices, $v\equiv 1\pmod{4}$, then $D(G)$ has a signature $( \frac{v+1}{2}, \frac{v-1}{2},0)$. So for $v\ge 9$,  $n_{+}(G) > n_{-}(G)$.

In this note, we determine the distance signature for the incidence graph of an affine design.

\begin{thm}\label{th:main}
Let $\sD$ be an affine $(v,k, \lambda)$-design with $v$ points and $b$ blocks, where $v=n^2\mu,$ $b=(n^3\mu-1)/(n-1)$ for integers $n\ge 2$, $\mu\ge 1$. If $G$ is the incidence graph  of $\sD$, then  the distance signature  of $G$ is 
$$
\sig (G) = \begin{cases}
(1,4,5) & \mbox{ if } (n,\mu)=(2,1)\\
(4\mu,4\mu-1, 4\mu-2) &\mbox{ if } n =2, \mu \ge 2\\
(b,v,0) &\mbox{ otherwise.}
\end{cases}
$$
\end{thm}

Recently, 
 Zhang \cite{Zhang} determined the distance signatures of complete $k$-partite graphs, and Zhang and Godsil \cite{Zhang13} gave some graphs with $n_{+}(G)=1$ and obtained their distance signatures.  See the recent survey paper \cite{Aouc} for more background and activities about the spectra of distance matrices.

\section{Preliminaries}
A design $\mathscr{D}$ is a pair $(V, \mathcal{B})$, where $V$ is a finite set and and $\mathcal{B}$ is a set of subsets of $V$. Elements of $V$ and $\mathcal{B}$ are called points and blocks, respectively. $\sD$ is called a $(v, k, \lambda)$ design if (1) $|V|=v$, (2) every block contains $k$ points,  (3) every 2-subset of $V$ is contained in precisely $\lambda$ blocks. For a $(v, k, \lambda)$ design, it can be shown that the number of blocks containing any fixed point is $\lambda(v-1)/(k-1)$, denoted by $r$. The number of blocks $b$ is $vr/k$. People sometimes refer $\sD$  to as  $(v,b,r,k,\lambda)$ design.

A {\em parallel class} in $\mathscr{D}$ is  a set of blocks that form a partition of  $V.$  $\sD$ is called a {\em resolvable} design if $\mathcal{B}$ admits a partition of parallel classes:
$$
\mathcal{B} = \Cc_1 \cup \Cc_2\cup \dots \cup \Cc_r.
$$
Blocks from the same parallel class are said to be {parallel}. 
A resolvable $(v,k,\lambda)$ design $\sD$ is called  {\em affine }  if 
any two distinct non-parallel blocks intersect at exactly $\mu$ points. Each parallel class has $n=v/k$ blocks. The parameters $v,b,k,r,\lambda$ can be expressed in terms of $n$ and $\mu$:
\begin{equation} \label{eq:ad}
v = nk = n^2\mu, \quad k = n\mu, \quad 
\lambda = \frac{n\mu -1}{n-1}, \quad r = \frac{n^2\mu -1}{n-1}, \quad 
b = v+r -1 = \frac{n(n^2 \mu -1)}{n-1}.
\end{equation}
We also denote this design by AD($n,\mu)$. In particular, if $\mu =1$, $\sD$ is called an {\em affine plane} of  order  $n$.  

Examples of affine designs are given by the affine spaces. Let $V$ be an $m$-dimensional vector space over a finite field $\mathbb{F}_q$. A coset of an $i$-dimensional subspace is called $i$-flat, and the 0-, 1-, 2-, and $(n-1)$-flats are called points, lines, planes, and hyperplanes, respectively.  Now the points and hyperplanes form an affine $(q^m, q^{m-1}, \frac{q^{m-1}-1}{q-1})$ design, denoted by AG$(m,q)$, with $n=q$,  $\mu=q^{m-2}$. (Blocks are parallel if and only if they are cosets of the same subspace. If $B, B'$ are non-parallel blocks, then $B\cap B'$ is a coset of an $(m-2)$-dimensional subspace, and so has cardinality $q^{m-2}$.) 

  AG$(2,q)$ is  an affine plane and  AG$(m,2)$ is an AD$(2,2^{m-2})$.
Further examples are given by the Hadamard designs. Every Hadamard $(4\mu-1, 2\mu-1,\mu-1)$ design has an affine extension AD$(2,\mu)$, in which each parallel class consists of a block and its complement. (This extension is called a Hadamard 3-design.)  Moreover, 
  an AD($2,\mu)$ is equivalent to a Hadamard $(4\mu-1, 2\mu-1,\mu-1)$ design or a Hadamard matrix of order  $4\mu$.  See  \cite[Ch. 1]{Cameron} or \cite[Sec. II.8]{Beth} for more details.   

Now we recall the incidence matrix of a design. In the sequel, we let $I_s$  be the identity matrix of order $s$, and    $J_{s,t}$ (or $O_{s,t}$) be the $s\times t$ matrix of all-one (or all-zero) entries.  We also write $J_s = J_{s,s}$ and  $O_s = O_{s,s}$.  The subscripts $s,t$  may be suppressed when they are clear from the context.
Let $\sD= (V, \mathcal{B})$ be a $(v,k,\lambda)$ design.
The {\em incidence matrix}\footnote{
Some authors define the incidence matrix as the transpose of $N$.}%
   $N$ of $\sD$ is formed by indexing the rows with points $V$ and columns with blocks $\mathcal{B}$ and defining the $(p,b)$ entry  to be $1$ if  $p\in b$ and $0$ otherwise. By definition, 
\begin{equation}
NN^t = r I_v + \lambda (J_v - I_v).
\end{equation}   
On the other hand,  $N^t N$ has $(b,b^\prime)$ entry $|b\cap b^\prime|$.

The {\em incidence graph} of design $\sD =(V, \mathcal{B)}$ is the bipartite graph on the vertex set $V\cup \mathcal{B}$ with independent sets $V$ and $\mathcal{B}$
and for $v\in V$, $b\in \mathcal{B}$, $v$ and $b$ are  adjacent if $v\in b$. So this graph has adjacency matrix 
$$
\begin{bmatrix} O & N \\ N^t & O\end{bmatrix}.
$$
For typographic convenience, we sometime leave blank a zero block.

\section{Proofs} \label{proofs}

In this section, we prove our main theorem. 

Suppose that  $\sD$ is an affine $(v,b,r,k,\lambda)$ design and these parameters are given as in Eq. \eqref{eq:ad}.
Let  $\Cc_1, \Cc_2, \dots, \Cc_r$, ($|\Cc_i|=n$), be a system of parallel classes   and $N$ be the incidence matrix of $\sD$. Then the product $N^t N$ is a $r\times r$ block matrix and each block is 
a square matrix of size $m$.
The diagonal blocks ($\Cc_i,\Cc_i)$ are $k I_n$ and the off-diagonal blocks ($\Cc_i,\Cc_j)$ are $\mu J_n$: 
\begin{equation}
N^t N = k I_r \otimes I_n + \mu (J_r - I_r) \otimes J_n,
\end{equation} 
where $A\otimes$ is the Kronecker product of matrices.  

The following equations hold by the definition of an affine design:
\begin{align}\label{eq:ad2}
J_v N = kJ_{v,b}, \quad  N^t J_v = k J_{b,v},\quad N J_b = r J_{v,b}, \quad J_b N^t = r J_{b,v}.
\end{align}

Let $G$ be the incidence graph of $\sD$ and $A_i$ be  the distance-$i$ matrix of $G$.
Now we have
$$
A_0 = I_{v+b}, \quad A_1 = A = \begin{bmatrix} O_v & N \\ N^t & O_b \end{bmatrix}.
$$
We calculate $A_i$ recursively from the products $AA_{i-1}$:
$$
\begin{array}{ll}
A^2   = \begin{bmatrix} NN^t & \\ & N^t{N} \end{bmatrix}
    =  \begin{bmatrix}
       rI_v + \lambda (J_v - I_v) &
       \\ &  k I_r \otimes I_n + \mu (J_r - I_r) \otimes J_n
         \end{bmatrix},
 \quad
 & A_2  = \begin{bmatrix}
          J_v - I_v & \\ & (J_r  - I_r)\otimes J_n \end{bmatrix}  ,
\\[15pt]
A A_2   = \begin{bmatrix}
          & N  (J_r  - I_r)\otimes J_n \\
            N^t (J_v - I_v) &
               \end{bmatrix}
               = \begin{bmatrix}
               & (r-1) J_{v,b} \\ k J_{b,v} - N^t &
             \end{bmatrix},
\quad
& A_3  =  \begin{bmatrix}
         &  J_{v,b} - N & \\  J_{b,v}  - N^t & \end{bmatrix}  ,
\\[15pt]
AA_3    
       = \begin{bmatrix}
       (r-\lambda)(J_v - I_v) & \\ & k I_r \otimes(J_n - I_n) + (k-\mu)(J_r -I_r)\otimes J_n
       \end{bmatrix},
\quad
& A_4 = \begin{bmatrix} O_{v} & \\ & I_r\otimes (J_n - I_n)
       \end{bmatrix}    .
\end{array}
$$
So graph $G$ has diameter $4$ with distance matrices $A_i \, (0\le i\le 4)$ given above.

Let $D = A_1 + 2A_2 + 3A_3 +4A_4$ be the distance matrix of $G$:
\begin{equation}\label{dM}
D = \begin{bmatrix}
   2(J_v - I_v) & 3J_{v,b} -2N
   \\
   3J_{b,v} - 2N^t & 2(J_r - I_r)\otimes J_n + 4 I_r \otimes (J_n -I_n)
\end{bmatrix}.
\end{equation}

Now we will determine the characteristic polynomial of $D$: 
$$
xI_{v+b} -D = \begin{bmatrix}
   x I_v - 2(J_v - I_v) & -(3J_{v,b} -2N)
   \\
   -(3J_{b,v} - 2N^t)  & xI_b - [2(J_r - I_r)\otimes J_n + 4 I_r \otimes (J_n -I_n)]
\end{bmatrix}
:= \begin{bmatrix}
C_0 & C_1 \\ C_1^t & C_2
\end{bmatrix}.
$$
Treating $x$ as indeterminate, $C_0 = (x +2)I_{v} - 2J_v  $ is invertible:
\begin{align*}
C_0^{-1} =
\dfrac{1}{(x+2)(x-2v+2)}
\begin{bmatrix}
x-2v+4 & 2 & \dots &2
\\
2 & x-2v+4& \dots & 2
\\
\vdots & \vdots& \ddots & \vdots
\\
2 & 2 & \dots & x-2v+4
\end{bmatrix}
= \dfrac{(x-2v+2)I_v + 2J_v}{(x+2)(x-2v+2)}.
\end{align*}
and $C_0$ has determinant
\begin{equation}\label{eq:C0}
\det C_0 = (x-2v+2)(x+2)^{v-1}.
\end{equation}

Since 
$$
\begin{bmatrix}
I_b & \\ C_1^t C_0^{-1}  & I_v
\end{bmatrix}
\begin{bmatrix}
C_0 & C_1 \\ C_1^t & C_2
\end{bmatrix} 
= \begin{bmatrix}
C_0 &  \\[5pt]
  & C_2 - C_1^t C_0^{-1} C_1
\end{bmatrix} ,
$$
we have 
\begin{equation}\label{eq:det}
\det \begin{bmatrix}
C_0 & C_1 \\ C_1^t & C_2
\end{bmatrix}
= \det C_0 \cdot \det ( C_2 - C_1^t C_0^{-1} C_1).
\end{equation}

Now from Eq. \eqref{eq:ad}-\eqref{eq:ad2} we obtain 
\begin{align*}
C_1^t C_0^{-1} C_1 &= \dfrac{1}{(x+2)(x-2v+2))}
(2N^t - 3J_{b,v})
[(x-2v+2)I_v + 2J_v) ]
(2N - 3J_{v,b})
\\ 
& = \dfrac{1}{(x+2)(x-2v+2)}
[4(x-2v+2)N^t N + J_b(8k^2 +  9vx + 18v -12kx -24k)],
\end{align*}
\begin{align*}
C_2 - C_1^t C_0^{-1} C_1 =&
I_r\otimes I_n \cdot \left(x- \frac{4k}{x+2} -
\frac{8k^2 +9vx +18v-12kx -24k}{(x+2)(x-2v+2)}\right)
\\ %
 &- I_r\otimes (J_n -I_n)\cdot \left( 4 +
 \frac{8k^2 +9vx +18v-12kx -24k}{(x+2)(x-2v+2)}\right)
\\ %
 & -(J_r-I_r)\otimes J_n \cdot \left( 2 + \frac{4\mu}{x+2} +
  \frac{8k^2 +9vx +18v-12kx -24k}{(x+2)(x-2v+2)}\right)
 \\
 :=&\ \alpha I_r \otimes I_n + \beta I_r \otimes (J_n - I_n) + \gamma (J_r - I_r)\otimes J_n.
\end{align*}
 Note that  $C_2 - C_1^t C_0^{-1} C_1$  is a $r\times r$ block matrix, with respective diagonal and off-diagonal blocks
$$
\begin{bmatrix}
 \alpha &\beta&\dots &\beta\\
   \beta&\alpha &\dots &\beta\\
   \vdots &\vdots &\ddots &\vdots\\
     \beta&\beta&\dots &\alpha
\end{bmatrix}
\quad\mbox{and}
\begin{bmatrix}
 \gamma& \gamma& \dots &\gamma \\
  \gamma& \gamma& \dots & \gamma\\
  \vdots &\vdots &\ddots &\vdots\\
  \gamma& \gamma& \dots & \gamma
\end{bmatrix}.
$$

\begin{lemma}\label{lem:D}
Let  $D(r,m,\alpha,\beta,\gamma) = \det(C_2 - C_1^t C_0^{-1} C_1)$. Then we have
 $$
 D(r,m,\alpha ,\beta ,\gamma) =
  (\alpha +(m-1)\beta  +m(r-1)\gamma)
  \cdot  (\alpha +(m-1)\beta  - m\gamma)^{r-1}
  \cdot
 (\alpha -\beta )^{r(m-1)} .
 $$
\end{lemma}
\begin{pf}
Now we perform the following operations on  $D(r,m,\alpha,\beta,\gamma)$:
First adding all rows to the first row, the first row has the same entry
$\alpha + (m-1)\beta + m(r-1)\gamma$. Now factoring this number out so the resulted matrix has all 1 in the first row.  Add a suitable multiple of the first row to each of rest rows to create following determinant (eg. $r=3$ below):
$$
\left|\begin{array}{cccc|cccc|cccc}
 1  &1 &\dots &1 &  1  &1 &\dots &1 &  1  &1 &\dots &1 \\
   \beta&\alpha &\dots &\beta &   \gamma& \gamma& \dots & \gamma &  \gamma& \gamma& \dots & \gamma\\
   \vdots &\vdots &\ddots &\vdots &   \vdots &\vdots &\ddots &\vdots &   \vdots &\vdots &\ddots &\vdots\\
     \beta&\beta&\dots &\alpha &   \gamma& \gamma& \dots & \gamma  &   \gamma& \gamma& \dots & \gamma
\\\hline
 0&0& \dots &0 &  \alpha - \gamma &\beta- \gamma &\dots &\beta - \gamma & 0&0& \dots &0 \\
 0& 0 & \dots &0 &   \beta- \gamma&\alpha - \gamma &\dots &\beta- \gamma &  0& 0 & \dots &0\\
  \vdots &\vdots &\ddots &\vdots & \vdots &\vdots &\ddots &\vdots &   \vdots &\vdots &\ddots &\vdots\\
  0 & 0& \dots & 0   &  \beta- \gamma&\beta- \gamma&\dots &\alpha- \gamma &    0 & 0& \dots & 0
\\\hline
 0&0& \dots &0 &0&0& \dots &0 &  \alpha - \gamma &\beta- \gamma &\dots &\beta - \gamma\\
 0& 0 & \dots &0 & 0& 0 & \dots &0 &  \beta- \gamma&\alpha - \gamma &\dots &\beta- \gamma\\
  \vdots &\vdots &\ddots &\vdots &  \vdots &\vdots &\ddots &\vdots&    \vdots &\vdots &\ddots &\vdots\\
  0 & 0& \dots & 0    &   0 & 0& \dots & 0       &  \beta- \gamma&\beta- \gamma&\dots &\alpha- \gamma
\end{array}\right|.
$$
The first diagonal block has determinant $(\alpha - \beta)^{m-1}$ and the second diagonal block has determinant
$ \alpha + (m-1)\beta - m\gamma$. Now the proof is completed. \hfill $\Box$
\end{pf}

From Eq. \eqref{eq:ad}, we find
\begin{align}
\alpha - \beta = x +4 -\frac{4n\mu}{x+2}, \quad \alpha +(m-1)\beta  - m\gamma = x-2n+4,\label{eq:ab1}
\\ \notag
\alpha +(m-1)\beta  +m(r-1)\gamma=
(x-2v+2)^{-1}\left( x^2  -\frac {2x(2\,
{n}^{3}\mu-\mu\,{n}^{2}+{n}^{2}-5\,n+3)}{n-1}\right.
\\ \left. 
-{\frac { \left( \mu\,{n}^{2}-1 \right)  \left( 5\,{n}^{3}\mu-12\,\mu
\,{n}^{2}+4\,\mu\,n-4\,{n}^{2}+16\,n-8 \right) }{n-1}}
\right)  \label{eq:ab2}
\end{align}
and
\begin{equation}\label{eq:rel}
r(n-1)  = v-1, \quad 
\end{equation}

The characteristic polynomial of $D$ follows from Lemma \ref{lem:D} and Eq. \eqref{eq:C0}-\eqref{eq:rel}.
\begin{thm}\label{th:poly}
Let $D$ be the distance matrix in \eqref{dM}. Then the characteristic polynomial of $D$ is 
\begin{align*}
(x-2n+4)^{r-1} (x^2+6x+8-4n\mu)^{v-1} \left( x^2  -\frac {2x(2\,
{n}^{3}\mu-\mu\,{n}^{2}+{n}^{2}-5\,n+3)}{n-1}\right.
\\ \left. %
-{\frac { \left( \mu\,{n}^{2}-1 \right)  \left( 5\,{n}^{3}\mu-12\,\mu
\,{n}^{2}+4\,\mu\,n-4\,{n}^{2}+16\,n-8 \right) }{n-1}}
\right).
\end{align*}
\end{thm}
\vspace{10pt}

\textbf{Proof of Theorem \ref{th:main}.}
Consider roots for each factor in the above theorem, the superscripts in parentheses indicating multiplicities: 
 $(x-2n+4)^{r-1}$ has roots $0^{(r-1)}$ if $n=2$, and only positive root $(2n-4)^{(r-1)}$ otherwise;
 $(x^2+6x+8-4n\mu)^{v-1}$ has roots $0^{(v-1)}, -6^{(v-1)}$ if $(n,\mu)=(2,1)$, and $-3 \pm \sqrt{1+4n\mu}^{(v-1)}$ otherwise;
 it is a bit tedious, but still easy to see that the last factor has one positive and one negative root for any $n\ge 2$, $\mu\ge 1$. Now Theorem \ref{th:main} follows. \hfill $\Box$

\section*{Acknowledgments}
The author is supported in part  by the Science Foundation of  Hebei Normal University (grant No. L2011B02), and by the National Natural Science Foundation of China  (grant No. 11571091).

This work was started  while I visited Professor Eiichi Bannai at Shanghai Jiaotong University from Sept. 25 to Oct. 3, 2014. He mentioned the conjecture of Dr. Kohei Yamada about the distance signatures of affine designs. I solved the problem for  affine planes while there and  Theorem \ref{th:poly}  one week later, but I did not type up the details.   When I emailed Bannai my solution in mid January 2015, we learned that Yamada also proved his conjecture. His statement of Theorem \ref{th:poly} is of a different form. In fact, Yamada had determined the distance signatures for many other designs as well.  
 I would like to thank Professor Bannai for many communications, Dr. Yamada for sharing his results, and 
 the anonymous referee for a helpful report. 
 
 
\end{document}